# I. Sh. Jabbarov

## On mean values and zeroes of Dirichlet series

### 1. Introduction

In this paper the mean values and zeroes of Dirichlet series of a view

$$f(s) = \sum_{n=1}^{\infty} a_n n^{-s}, \qquad (1)$$

with complex coefficients $a_n$ are studying; here $s = \sigma + it$ is a complex variable. It is known that if the series (1) is convrgent at the point $s = \sigma_0 + it$ then in the half plane $\sigma > \sigma_0$ the function $f(s)$ is regular (see [22, p.326]). Most of analytical properties of Dirichlet series are connected with mean values of Dirichlet series (see [4, 6, 10, 18, 21, 22, 23, 25, 26, 28, 34, 37, 40]). If there exist a finite limsup

$$M(f) = M(f,\sigma) = \varlimsup_{T \to \infty} \frac{1}{T} \int_0^T |f(\sigma + it)|^2 dt \quad (2)$$

then one calls it to be a mean value of the Dirichlet series (1). The basic question consists in finding of half planes $\sigma > \alpha$, with possible smaller real $\alpha$, on which (2) exists. It is known, for example, that the limsup (2) exists in the half plane of absolute convergence (see [22, p.340]). Consideration of mean values for Dirichlet series in wider areas where to the series (1) can be analytically continued, generally, is connected with great difficulties.

In [22, chapter IX] mean values of series having a finite order are considered. For such Dirichlet series the existence of (2) involves the convergence of the series (3) (see below). In the same place, the function $\mu(\sigma)$ designating an order of Dirichlet series on the line $\mathrm{Re}\, s = \sigma$ is entered [22, p. 335]. Thus, the notion of half plane of mean values associates with an existence of finite order on this half plane. Well known Lindelöf Hypothesis for the Riemann zeta-function asserts that $\mu(\sigma) = 0$ in the half plane of mean values (see [22, p. 328]).

Only in some cases an existence of limsup (2) outside the half plane of absolute convergence was obtained. It is known, for example, that at $\sigma > 1/2$

$$\lim_{T \to \infty} \frac{1}{T} \int_0^T |\zeta(\sigma + it)|^2 dt = \zeta(2\sigma)$$

(see [23]). The more stronger result of Hardy and Littlewood (see [23, p. 139])

$$\lim_{T \to \infty} \frac{1}{T} \int_0^T |\zeta(\sigma + it)|^4 dt = \frac{\zeta^4(2\sigma)}{\zeta(4\sigma)},$$

held in $\sigma > 1/2$, shows an existence of mean values of the square of Riemann zeta-function in the specified half plane. For higher degrees of the zeta-function the area of existence of limsup (2) is narrow (see [23, chap. VII]). Similar results are known for Dirichlet $L$-functions.

The Lindelöf Hypothesis is equivalent to the statement

$$\lim_{T \to \infty} \frac{1}{T} \int_0^T |\zeta(\sigma + it)|^{2k} dt = \sum_{n=1}^{\infty} \tau_k^2(n) n^{-2\sigma},$$

for any natural $k$ and $\sigma > 1/2$ (see [23, p. 328]). On the left part of the last equality the mean values of $k$-th degree of the zeta-function stands.

In the book [22] the communication of mean values of Dirichlet series with the convergence is studied. The theorem asserting is proved that the Dirichlet series (1) converges in the half plane of mean values, if it represents a regular function on this half plane.

In the present work somewhat the converse problem in more general conditions is considered. Namely, considering the half plane of mean values, we define it as the half plane $\sigma > \sigma_m$ where the series

$$\sum_{n=1}^{\infty} |a_n|^2 n^{-2\sigma} \quad (3)$$

converges, do not assuming that the Dirichlet series defines a function of a finite order (see [19, p. 341]). We prove that for one class of Dirichlet series of a kind (1) from the convergence of the Dirichlet series in the above defined half plane it follows that it has there, not only a finite order (see [22, p. 335]), but there is the mean value (2) at $\sigma > \sigma_m$ for any natural degree of the function $f(s)$ also. At first, we define a class of Dirichlet series satisfying (3) for which we will establish our results. Let $r \geq 2$ be a natural number and $N(r)$ designates the set of all those natural numbers the canonical factorizations of which contain only the prime numbers $p \leq r$. We say that the Dirichlet series (1) belongs to the class $J = J(\sigma_m)$ if the following series converges absolutely at $\sigma > \sigma_m$:

$$\sum_{n \in N(r)} a_n n^{-s}.$$

Let's notice that the majority of often used Dirichlet series as the zeta-function or Dirichlet $L$-functions belongs to this class.

**Theorem 1.** *Let the function $f(s)$ being defined by the series (1) belong to the class*

$J = J(\sigma_m)$, *and be convergent in the half plane* $\sigma > \sigma_m$ *where the series below converges:*

$$\sum_{n=1}^{\infty} |a_n|^2 n^{-2\sigma}.$$

*Then in the half plane* $\sigma > \sigma_m$:

*1) for any natural k the following relation holds:*

$$\lim_{T \to \infty} \frac{1}{T} \int_0^T |f(\sigma + it)|^{2k} dt < +\infty;$$

*2)* $\mu(\sigma) = 0$.

The distribution of zeroes of a given Dirichlet series is a problem the deep investigation and getting of expected results of which did not succumb to the classical methods of analysis. Easy observation shows that every Dirichlet series has a half plane, placed in the half plane of absolute convergence, free from the zeroes (see [22, p.347]). Further, if the Dirichlet series has mean values then, denoting by $N(\sigma,t)$ the number of zeroes $\sigma' + it'$ of $f(s)$ satisfying the conditions $\sigma' > \sigma, |t| \leq T$, we have $N(\sigma,T) = O(T)$ (see [22, p.349]).

As it is clear from the discussions of n° 9.7 of the book [22, p. 351], the Dirichlet series has in the half plane of absolute convergence a property of almost periodicity. The results of Voronin S. M. show that in the half plane of mean values non-zero values of Dirichlet series (see[25]) are distributed, likely as values of the almost periodic functions are distributed. But the method of Voronin S. M. is non-successful in the question on the distribution of zeroes. Montgomery H. L. in [18, 19] had shown that every zero of the zeta function induces a number of additional zeroes placed in some rectangle containing the given zero.

Considering the class $J(\sigma_m)$ of Dirichlet series, we show that the series of this class have not zeroes in the half plain of mean values, if it is regular in this half plain and have inverse Dirichlet series belonging to the $J(\sigma_m)$. For this we show that, firstly,

$$N(\sigma,T) = o(T), \quad (4)$$

and, secondly, every possible zero in the half plane of mean values induces more than $\gg T$ zeroes on the same line (in the considered rectangle). The results of the first type are known as the density theorems on the distribution of zeroes of Dirichlet series. More information about history and results of this theory can be found in [4,6,7,10,17,18,21,23,33]). First result in this direction was received in [41] where the estimation $N(\sigma,T) = O(T)$ was proven for the zeta function. This result of Bohr H. and Landau E. (wich later was proved in the form (4)) for the zeta function was improved by Carlson F. and further investigations on the same direction were performed after of discovering by Hoheisel G. the possibility of the use of density theorems on the zeroes' distribution for arithmetical applications. Yu. V. Linnik ([17]) gave a new proof of the theorem of Holgbach - Vinogradov on the representation of great even numbers as a sum of three primes by the method of Hardy and Littlewood using density theorems on zeroes of the zeta function as a substitute for the Riemann Hypothesis. One of the most important consequences of this theory is a theorem of Bombieri E. and Vinogradov A. I. which successfully used in many questions of the theory of prime numbers.

The first result of the second type at first was gotten by Voronin S. M. in [29] for the zeta functions of quadratic forms. The basic moment here is that fact that this function can be represented as a linear combination of some Dirichlet series with Euler product. The similar method was applied in [25] for investigation of non zero values of the zeta function. We show that the method developed by Voronin S. M. is possible to extend for investigation of the zeroes of Dirichlet series belonging to the class $J(\sigma_m)$ eliminating aroused difficulties. For us the weak result (4) indicated above is sufficient. The main result looks as the theorem below.

**Theorem 2.** *Let the function $f(s)$, being defined by the series (1), belong to the class $J(\sigma_m)$ and be convergent in the half plane $\sigma > \sigma_m$. Let in the half plane $\sigma > \bar{\sigma}$ of absolute convergence there exist $f(s)^{-1} = \sum_{n=1}^{\infty} b_n n^{-s}$ in some half plane contained in the half plane of absolute convergence and $f(s)^{-1} \in J(\sigma_m)$. Then $f(s) \neq 0$ in the half plane $\sigma > \sigma_m$.*

The theorems' proofs based on some new results in the metric theory of distribution of curves in the infinite dimensional unite cube $\Omega = [0,1] \times [0,1] \times \cdots$ where a new measure (see[11-14]) is constructed.

## 2. Introdusing of a new measure

It is best known that in the cube $\Omega$ a product Lebesgue measure may be introduced (see [5, p. 219]). There is, also, another construction of a measure in $\Omega$ called the Haar measure. The Haar measure is a measure defined in the locally compact topological groups. It was proven also uniqueness of this measure (see [1, p. 303], [8, p.241]). Many of measures used in various brunches of the mathematics could be considered as a Haar measure. Particularly, the product of Lebesgue measures in $[0,1]^n$ is a Haar measure for any natural number $n$ and, hence, is unique in this cube. Really, to prove this, consider the topological group $R^n$. The group $Z \oplus \cdots \oplus Z$ is a subgroup, so the factor group $T_n = R^n /(Z \oplus \cdots \oplus Z)$ as a compact group is locally compact. Therefore, the invariant measure in this group is unique. Let $A \subset [0,1]^n$ is a measurable set in the product Lebesgue meaning in $[0,1]^n$. Consider the union of intersections $(\bar{a} + A) \cap (\bar{m} + [0,1)^n)$, $\bar{m} \in Z^n$ for any given vector $\bar{a} \in R^n$. Only no more than $2^n$ of these intersections are non-empty, and the sum of their measures is equal to the measure of the set $A$. Therefore, the product measure is invariant in regard to the transitions $\bar{x} \mapsto \bar{x} + \bar{a} \pmod 1$, $\bar{x} \in [0,1)^n, \bar{a} \in R^n$. So, this is a unique measure with the property of invariance.

Despite that that the most of told above is true for $\Omega$, the situation is currently different in the infinite dimensional case. We can define the Haar measure in $\Omega$, as in the factor group, invariant in regard to the transitions (mod 1). We get, then, some unique measure defined in $R^\infty$.

Consider now the unite cube $\Omega$ itself. In this case the "number" of non-empty intersections $(\bar{a} + A) \cap (\bar{m} + \Omega')$, $A \subset \Omega'$ (here $\Omega' = \{(\omega_n) \mid 0 \leq \omega_n < 1\}$) is non-countable. So, we can not state that the product measure in $\Omega$ is invariant in regard to the transitions $\bar{x} \mapsto \bar{x} + \bar{a} \pmod 1, \bar{x} \in \Omega', \bar{a} \in \Omega$. Therefore, in $\Omega$ it could be introduced the measure different from the Haar measure restriction of which in $\Omega$ coinsides, partially, with the product Lebesgue measure. Some of sets being measurable in the Haar meaning can stand now nonmeasurable. Morover, another measure, different from the product Lebesgue measure in $\Omega$ could be introduced also. The lemma below shows justness of these statements. To formulate it we need in some designations.

We begin with introducing of special curves of a kind $(\{t\lambda_n\})_{n \geq 1}$ (the sign $\{\}$ means a fractional part, and $\lambda_n > 0, \lambda_n \to \infty$ as $n \to \infty$) in the subsets of infinite dimensional unite cube. In the works [2, 15, 20, 22] the finite case has been studied. It is necessarily to note that this curve has a zero measure in the product Lebsgue or Haar meaning. The basic result of the present section is formulated in the lemma 1 below which based on the lemma 2 used in particular cases in the works ([11-13]). We recall some designations and definitions of those works.

**Definition 1.** *Let $\sigma : N \to N$ be any one to one mapping of the set of natural numbers. If for any $n > m$ there is a natural number $m$ such that $\tau(n) = n$, then we call $\tau$ to be a finite permutation. A subset $A \subset \Omega$ is called to be finite-symmetrical if for any element $\theta = (\theta_n) \in A$ and any finite permutation $\tau$ one has $\tau\theta = (\theta_{\tau(n)}) \in A$.*

Let $\Sigma$ to denote the set of all finite permutations. $\Sigma$ is a group which contains each group of $n$ degree permutations as a subgroup (we consider each $n$ degree permutation $\tau$ as a finite permutation, in the sense of definition 1, for which $\tau(m) = m$ when $m > n$). The set $\Sigma$ is a countable set and we can arrange its elements in a sequence.

Let $\omega \in \Omega$, $\Sigma(\omega) = \{\tau\omega \mid \tau \in \Sigma\}$ and $\Sigma'(\omega)$ means the closed set of all limit points of the sequence $\Sigma(\omega)$. For the real $t$ we denote $\{t\Lambda\} = (\{t\lambda_n\})$ where $\Lambda = (\lambda_n)$. Let $\mu$ to denote the product of linear Lebesgue measures of $m$ given on the interval $[0,1]$: $\mu = m \times m \times \cdots$.

In $\Omega$ we define the Tychonoff metric by the following expression

$$d(x,y) = \sum_{n=1}^{\infty} e^{1-n} |x_n - y_n|.$$

Then, $(\Omega, d)$ becomes a compact metric space. The $\mu_0$- measure of the ball was defined in [12,13,14]. On this bases it may be introduced a measure $\mu_0$ in the $\Omega$ by known way by using of open sets: the open ball we define as an intersection $\Omega \cap B(\theta, r)$; an elementary set we define as a set being gotten by using of finite number of operations of unionize, taking differences or complements. It is clear that every elementary set is $\mu_0$-measurable. The set of elementary subsets of $\Omega$ is an algebra of subsets. The $\sigma$-algebra of subsets in $\Omega$ can be defined by known way and the set function defined in the subsets' algebra can be extended to the $\sigma$-algebra (see [5, p. 152]). The outer and inner measures can be introduced by known way. Given any set in $\Omega$, we call it to be measurable, if and only if, when it's outer and inner measures are equal (see also [1,VI, p.16]). Defined measure will be, as it seen from the reasoning above, a regular measure, and outer measure of a set in the meanings of product and $\mu_0$-meaning are the same. Really, every open set can be represented as a countable union of open balls. Since, as it was shown in [12,13], the measure of every open ball could be approximated with any small error by the volumes $(\mu_N(r))$ with suitable $N$. So, every measurable set, in the meaning of introduced measure, is measurable in the meaning of product Lebesgue measure also.

The projection of the curve $(\{t\lambda_n\})_{n\geq 1}$ in two dimensional plane, i. e. the curve $(\{t\lambda_1\}, \{t\lambda_2\})$ has a zero measure. By the theorem of Fubini, then, the product Lebesgue measure of the curve $(\{t\lambda_n\})_{n\geq 1}$ is also equal to zero.

**Lemma 1.** *Let the sequence $(\lambda_n)$ be an unbounded, monotonically increasing sequence of positive real numbers every finite subfamily of elements of which is linearly independent over the field of rational numbers. Then the curve $(\{t\lambda_n\}), t \in [0,1]$ is not $\mu_0$-measurable set in $\Omega$.*

The proof of this lemma is based on the following basic lemma of the works [12,13].

**Lemma 2**. *Let $A \subset \Omega$ be a finite-symmetric subset of zero measure and $\Lambda = (\lambda_n)$ is an unbounded, monotonically increasing sequence of positive real numbers any finite subfamily of elements of which is linearly independent over the field of rational numbers. Let $B \supset A$ be any open, in the Tychonoff metric, subset with $\mu_0(B) < \varepsilon$,*

$$E_0 = \{0 \leq t \leq 1 \,|\, \{t\Lambda\} \in A \wedge \Sigma'\{t\Lambda\} \subset B\}.$$

*Then, we have $m(E_0) \leq c_0 \varepsilon$ where $c_0 > 0$ is an absolute constant, $m$ designates the Lebesgue measure.*

*Proof of Lemma 1.* Let us suppose the converse statement. Let the curve $(\{t\lambda_n\}), t \in [0,1]$ be measurable. Then it measurable in the product meaning and has a zero measure. Therefore, the union $U = \bigcup_{0 \le t \le 1} \Sigma(\{t\Lambda\})$ as a set constructed from the curve $(\{t\lambda_n\}), t \in [0,1]$ by an action of the group $\Sigma$ of all finite permutations has a zero measure also, since it is a countable union of sets of zero measure. The set $U$ is finite-symmetrical. Let $n$ be any natural number. If we take a projection of the set $U$ into $\Omega$ by omitting the first $n$ coordinates (restricting of the sequence $(\{t\lambda_n\})$), we get again the set $U_n$ of zero measure. Really, the set $U$ can be overlapped by the union of balls with the total measure, not exceeding $\varepsilon > 0$. Restricting the ball $B(\theta_0, \lambda)$ by omitting the first $n$ coordinates, and denoting the projection by $S_N$, we get

$$S_N = \left\{ (\theta_n) \mid \sum_{n=N+1}^{\infty} |\theta_n - \theta_n^0| e^{1-n} < \lambda \right\}.$$

Since

$$\sum_{n=N+1}^{\infty} |\theta_n - \theta_n^0| e^{1-N} = e^{-N} \sum_{n=1}^{\infty} |\theta_{n+N} - \theta_{n+N}^0| e^{1-n},$$

then denoting the projection of the point $\theta_0$ by $\theta_0'$, we have $S_N = B(\theta_0', e^N \lambda)$, and $e^N \lambda \to 0$ as $\lambda \to 0$ for any fixed $N$. From this one deduces the demanded statement.

Consider the sequence of subsets $V_n = [0,1]^n \times U_n$ for all natural $n$. It is obviously that $V_n \subset V_{n+1}$. Let $V = \bigcup_{n=1}^{\infty} V_n$. We have $\mu(V_n) = 0$ for all values of $n$. Therefore, $\mu(V) = 0$ also, and the set $V$ is finite symmetrical. Then, there will be found some enumarable family of balls $B_r$ with the total measure, not exseeding $\varepsilon$, the union of which contains the set $V$. For every fixed natural $n$ we define the set $\Sigma_n'(t\Lambda)$ as a closed set of all limit points of the sequence (see [13]) $\Sigma_n(\overline{\omega}) = \{\tau\overline{\omega} \mid \tau \in \Sigma \wedge \tau(1) = 1 \wedge \cdots \wedge \tau(n) = n\}$ (i. e. the permutation $\tau$ remains the first $n$ indexes unchanged). Let

$$B^{(n)} = \{t \mid \{t\Lambda\} \in V \wedge \Sigma_n'(\{t\Lambda\}) \subset \bigcup_{r=1}^{\infty} B_r\}, \quad n = 1,2,....$$

For every $t$ the sequence $\Sigma_{n+1}(\{t\Lambda\})$ is a subsequence of the sequence $\Sigma_n(\{t\Lambda\})$. Therefore, $\Sigma_{n+1}'(\{t\Lambda\}) \subset \Sigma_n'(\{t\Lambda\})$ and we have $B^{(n)} \subset B^{(n+1)}$. Then, one gets the inequality $m(B) \le \sup_n m(B^{(n)})$ denoting $B = \bigcup_n B^{(n)}$. Let's estimate $m(B^{(n)})$ (we use the designations of the work [13, p.21]). We have

$$\{t\Lambda\} \in [0,1]^n \times \{\{t\Lambda\}'\} \subset V.$$

Let $(\theta_1,...,\theta_n) \in [0,1]^n$ be any point. There exist a neighborhood $V' \subset [0,1]^n$ of this point such that $(\theta_1,...,\theta_n,\{t\Lambda\}') \in V' \times W \subset \bigcup_r B_r$, for some neighborhood $W$ of the point $\{t\Lambda\}'$ (for this we can take, for example, the ball of half radiuse with the senter $\{t\Lambda\}'$). We, therefore, supplied every point $(\theta_1,...,\theta_n) \in [0,1]^n$ with some pair of open sets $(V',W)$. Since the set $[0,1]^n$ is closed, then it can be found a finite number of open sets $V'$ the union of which contains $[0,1]^n$. The intersection of corresponding open sets $W$, being an open set, contains the point $\{t\Lambda\}'$. Therefore, we have

$$[0,1]^n \times \{\{t\Lambda\}'\} \subset \bigcup V \times \bigcap W = [0,1]^n \times \bigcap W \subset \bigcup_{r \in R} B_r,$$

for each considered point $t$. The similar relationship is fair in the case when the point $\{t\Lambda\}$ would be replaced by any limit point $\overline{\omega}$ of the sequence $\Sigma(\{t\Lambda\})$ also, because $\overline{\omega} \in \bigcup_{r=1}^{\infty} B_r$. If one denotes by $B'$ the union of all open sets of the kind $\bigcap_{r \in R} B'_r$ (here $B'_r$ denotes the open set of trancated elements of $B_r$), corresponding to every possible values of $t$ and of the limit point $\overline{\omega}$, we shall receive the relation

$$\{t\Lambda\} \in [0,1]^n \times \{\{t\Lambda\}'\} \subset A \subset [0,1]^n \times B' \subset \bigcup_{r=1}^{\infty} B_r,$$

for each considered values of $t$ and

$$\{\overline{\omega}\} \in [0,1]^n \times \{\overline{\omega}\}' \subset A \subset [0,1]^n \times B' \subset \bigcup_{r=1}^{\infty} B_r,$$

for each limit point $\overline{\omega}$. From this it follows the inequality

$$\mu_0^*\left([0,1]^n \times B'\right) = \mu_0^*(B') \leq \varepsilon,$$

where $\mu_0^*$ means an outer measure. The set $B'$ is open and $\Sigma'(\{t\Lambda\}') \in B'$. Now we can apply the lemma 2 and receive the estimation $m(B^{(n)}) \leq c_0 \varepsilon$. Thus, we have $m(B) \leq c_0 \varepsilon$. Since $\varepsilon$ could be chosen arbitrarily small then there exist $t$ such that $t \notin B$. Then, $t \notin B^{(k)}$ for every $k = 1,2,...$. Consequently, for every $k$, there is a limit point $\overline{\omega}_k \in \Omega \setminus \bigcup_r B_r$ of the sequence $\Sigma_k(\{t\Lambda\})$. As the set $\Omega \setminus \bigcup_r B_r$ is closed, the limit point $\overline{\omega} = (\{t\Lambda\})$ of the sequence $(\overline{\omega}_k)$ will belong to the set $\Omega \setminus \bigcup_r B_r$. Therefore, $\{t\Lambda\} \notin \bigcup_{r \geq 1} B_r$ in the contradiction with our supposing. Then the curve $(\{t\lambda_n\})$, $t \in [0,1]$ could not be $\mu_0$-measurable. The proof of the lemma 1 is finished.

### 3. Basic auxiliary lemmas.

**Lemma 3.** *Let a series of analytical functions*

$$\sum_{n=1}^{\infty} f_n(s)$$

*be given in one-connected domain $G$ of a complex $s$-plane, and be absolutely converging almost*

*everywhere in G in Lebesgue sense, and the function*

$$\Phi(\sigma,t) = \sum_{n=1}^{\infty} |f_n(s)|$$

*is a summable function in G. Then the given series converges uniformly in any compact subdomain of G; in particular, the sum of this series is an analytical function in G.*

Proof of this lemma given in the works [14, 46].

Let $\alpha > \sigma_m, \delta = \alpha - \sigma_m, \sigma > \alpha$ and

$$f_k(s) = \sum_{n \in N(2^k)} a_n n^{-s},$$

where $k=1,2, \ldots$. Take the canonical factorization of the number $n$: $n = \prod_{p \backslash n} p^{\alpha_p}$. We have

$$f_k(s) = \sum_{n \in N(2^k)} a_n n^{-\sigma} e^{-2\pi i t \sum_{p \backslash n} v_p \alpha_p},$$

where $v_p = (\log p)/2\pi$.

**Lemma 4.** *Let $0 < r \leq \sigma - \alpha$ be any real number. For each real t there exists a limit*

$$\lim_{k \to \infty} f_k(s + \sigma + it) = f(s + \sigma + it)$$

*uniformly in the disc $|s| \leq r$.*

*Proof.* The elements of the infinite dimensional unite cube $\Omega = [0,1] \times [0,1] \times \cdots$ we will index by the prime numbers as $\theta = (\theta_p) \in \Omega$. Consider the trigonometrically series

$$g_k(\theta;s) = \sum_{n \in N(2^k)} a_n n^{-s} e^{-2\pi i \sum_{p \backslash n} \alpha_p \theta_p}; g_0(\theta;s) = a_1, |s| \leq r' = r + \delta/20.$$

From the uniqueness of canonical factorization one deduces for each complex $s$ that

$$\int_{\Omega} |g_k(\theta;s) - g_{k-1}(\theta;s)|^2 \mu(d\theta) = \sum_{n: p \wedge 2^{k-1} < p \leq 2^k} |a_n|^2 n^{-2\sigma} \leq \sum_{n=1}^{\infty} |a_n|^2 n^{-2\sigma}, \quad (5)$$

where the first summation is taken over all such $n$ every prime divisor of which does not exceed $2^k$ and has a prime divisor $p$ from the interval $2^{k-1} < p \leq 2^k$. Consider now the integral

$$h_k(\theta, \sigma) = \iint_{|s|=|\beta+i\tau| \leq r'} |g_k(\theta;s+\sigma) - g_{k-1}(\theta;s+\sigma)| d\beta d\tau, k \geq 1.$$

Let $C(\sigma) = \sum_{n=1}^{\infty} |a_n|^2 n^{-2\sigma}$ and $h(\theta,\sigma) = \sum_{k \geq 1} h_k(\theta,\sigma)$. In the conditions of the theorem we have $h(\theta,\sigma) \in L_1(\Omega, \mu_0)$. Really, from (5), applying Cauchy's inequality, we receive

$$\int_{\Omega} |h(\theta,\sigma)| \mu_0(d\theta) \leq r' \sqrt{\pi} \sum_{k=1}^{\infty} \left( \int_{\Omega} \iint_{|\beta+i\tau| \leq r'} |g_k(\theta; \beta+i\tau+\sigma) - g_{k-1}(\theta; \beta+i\tau+\sigma)|^2 d\beta d\tau \mu_0(d\theta) \right)^{1/2} \leq$$

$$\leq \pi r'^2 \sum_{k=1}^{\infty} \left( \sum_{n \in N(2^k) \setminus N(2^{k-1})} |a_n|^2 n^{-2\sigma+2r'} \right)^{1/2} \leq \pi r'^2 \left( 1 + \sum_{k=1}^{\infty} 2^{-(k-1)\delta/10} \right)^{1/2} C^{1/2}(\alpha - 3\delta/4). \qquad (6)$$

Therefore, from the theorem of Fatou (see [22, p. 387]) it follows that the series $h(\theta, \sigma)$ converges almost everywhere in the $\Omega$ and $h(\theta, \sigma) \in L_1(\Omega, \mu_0)$. From this, as a consequence, we deduce that the measure of a set $A_0$ of divergence of the series $h(\theta, \sigma)$ is equal to zero.

We can notice that if $\theta \in A_0$ then $\varphi\theta \in A_0$ for each finite permutation $\varphi \in \Sigma$, because of the divergence of permuted series. Now we follow by the construction used for the proof of the classical Egoroff's theorem (see [22, p. 379]). Let $(\varepsilon_j)$ be a sequence of positive numbers, $\varepsilon_j \to 0$. Consider in the set of convergence $\Omega_1 = \Omega \setminus A_0$ the sequence $q_n = \sum_{k>n} h_k(\theta; \sigma)$. Denote $S_{i,l} = \{\theta \in \Omega_1 \mid n \geq i \Rightarrow q_n < \varepsilon_l\}$. We have $S_{1,l} \subset S_{2,l} \subset \cdots$ and $\Omega_1 = \bigcup_{n \geq 1} S_{n,l}$ for each natural $l \geq 1$. Let $\Delta > 0$ be any positive number. One can define a natural number $n(l)$ for which $\mu_0(\Omega_1 \setminus S_{n(l),l}) < 2^{-l}\Delta$. Designate now $U_i = \bigcap_{l=i}^{\infty} S_{n(l),l}$. In every set $U_i$ the series $\sum_{k \geq 1} h_k(\theta; \sigma)$ converges uniformly and

$$\mu_0(\Omega_1 \setminus U_k) \leq \sum_{l \geq k} \mu_0(\Omega_1 \setminus S_{n(l),l}) < 2^{1-k}\Delta.$$

Further, noting that $U_1 \subset U_2 \subset \cdots$, we put $U = \bigcup_{k \geq 1} U_k$. So, $\mu_0(\Omega_1 \setminus U) = 0$ and if $\theta \in \Omega_1 \setminus U$ then $\theta \notin U_k$ for each $k$. Clearly, $A_0 \subset \Omega \setminus U$. Hence, for each natural $k$ one can find $l_k$ so that $\theta \notin S_{n(l_k),l_k}$. We can suppose that $n(l_k) \to \infty$ as $k \to \infty$, if else, then $n(l_k) \leq L$ for some natural $L$ and, so, $q_n < 2^{1-k}\delta$ for some $n$ and all $k$ which shows that $q_n = 0$. Then, the given series is absolutely convergent Dirichlet series (is a series of a view $f_r(s)$) for which the statement of the lemma 1 is true. So, we suppose that $n(l_k) \to \infty$. Now the relation $\theta \notin S_{n(l_k),l_k}$ shows that there exist a sequence of natural numbers $(r_k)$ such that $r_k \geq n(l_k)$, and $q_{r_k}(\theta, \sigma) \geq \varepsilon_{l_k}$. If $\varphi$ is a given finite permutation (which permutes only finite number of components $\theta_p$) then on great enough values of $k$ the components of $\theta$, taking part in the sum $q_{r_k}(\theta, \sigma)$, are invariable after of action of $\varphi$, so, $q_{r_k}(\varphi\theta, \sigma) \geq \varepsilon_{l_k}$. It means that $\varphi\theta \notin S_{n(l_k),l_k}$ for all great enough $k$, so that $\varphi\theta \notin U_{l_k}$. Since $l_k \to \infty$ then from this it follows that $\varphi\theta \notin U$, or, $\varphi\theta \in \Omega_1 \setminus U$, so, the last relation implies that the set $\Omega_0 = \Omega_1 \setminus U$ is finite symmetrical. *Moreover, if we should exchange the first $n(l_k) - 1$ components of the point $\theta$ by any numbers from the interval $[0,1]$ we get again some point*

$\theta' \in \Omega_0$. Since $n(l_k) \to \infty$ as $k \to \infty$, then this statement is fair for any first $n$ natural components.

Consider now the curve $(\{t\Lambda\})$ in $\Omega_0$. For every natural $n$ we define the set $\Sigma'_n(\omega), \omega \in \Omega$ as a set of all limit points of the sequence

$$\Sigma_n(\omega) = \{\sigma\omega \mid \sigma \in \Sigma \wedge \sigma(1) = 1, \wedge \cdots \wedge \sigma(n) = n\}.$$

Let

$$D^{(n)} = \{t \mid \{t\Lambda\} \in \Omega_0 \wedge \Sigma'_n(\{t\Lambda\}) \subset \bigcup_r B_r\}, \lambda_n = (1/2\pi)\log p_n, n = 1,2,\ldots;$$

here $p_n$ indicates the $n$-th prime number. Applying now the reasoning of the proof of the lemma 3 of [13], one deduces that the series

$$h(\theta,\sigma) = \sum_{k \geq 1} h_k(\theta,\sigma)$$

converges for almost all $\theta = \{t\Lambda\}, t \in \mathbf{R}$. Since such values of $t$ are everywhere dense in $\mathbf{R}$, the relation

$$\sum_{k=1}^{\infty} \iint_{|\beta+i\tau| \leq r+\delta/20} |g_k(\{t\Lambda\}; s+\sigma) - g_{k-1}(\{t\Lambda\}; s+\sigma)| d\beta d\tau = h(\{t\Lambda\},\sigma)$$

holds for all real $t$. Therefore, the series

$$\sum_{k=1}^{\infty} \iint_{|\beta+i\tau| \leq r+\delta/20} |f_k(s+\sigma+it) - f_{k-1}(s+\sigma+it)| d\beta d\tau \ ; f_0(s) = a_1$$

is convergent. From the theorem of Fatou it follows that the conditions of the lemma 3 are satisfied in the disc $|s| \leq r$ for the series $\sum_{k \geq 1}(f_k(\sigma+s+it) - f_{k-1}(\sigma+s+it))$ for each real $t$. Then, applying the lemma 3, we find

$$\lim_{k \to \infty} f_k(s+\sigma+it) = f(s+\sigma+it),$$

uniformly in s in the disc $|s| \leq r$ for all real $t$. Clearly, the speed of convergence is dependent on $t$. The lemma 4 is proven.

The following lemma is the best known Croneker's theorem (see [23, p.301], or [33, p.218]):

**Lemma 5.** *Let $\alpha_1, \alpha_2, \ldots, \alpha_N$ are real numbers linear independent over the field of rational numbers, $\gamma$ is a sub domain of $N$ –dimensional unite cube with the volume $\Gamma$ in Jordan meaning. Let further, $I_\gamma(T)$ is a measure of a set of such $t \in (0,T)$ for which $(\alpha_1 t, \alpha_2 t, \ldots, \alpha_N t) \in \gamma \pmod 1$. Then*

$$\lim_{T \to \infty} \frac{I_\gamma(T)}{T} = \Gamma.$$

The proof of this lemma can be found in [33, p. 345]. If the conditions of the lemma 5 is satisfied for any parallelepiped $\gamma$ then we say that the curve $\gamma(t) = (\gamma_n(t)) \in \Omega, t \geq 0$ uniformly distributed in $\gamma$ (mod 1) (see [33, p. 348]).

The following lemma is a generalization of the theorem of Croneker (see [33, p.345]).

**Lemma 6.** *Let the curve $(t\alpha_1, t\alpha_2, ..., t\alpha_N)$ be uniformly distributed (mod 1). Then for any integrable function $f(x)$, in Riemann sense, the following relation is true*

$$\lim_{T \to \infty} \frac{1}{T} \int_0^T f(\{t\alpha_1\}, ... \{t\alpha_N\}) dt = \int_0^1 \cdots \int_0^1 f(x_1, ..., x_N) dx_1 ... dx_N.$$

Under conditions of the lemma 6 the projection of the curve $(\{t\Lambda\})$ in the finite dimensional unite cube is uniformly distributed (mod 1) (see [33, p.348]).

**Lemma 7.** *For any continuous function $F(x)$ in the infinite-dimensional unite cube $\Omega$ and any open (or closed) ball $B$ the following relation*

$$\lim_{T \to \infty} \frac{1}{T} \int_{\{t\Lambda\} \in B} f(\{t\Lambda\}) dt = \int_B f(\theta) \mu(d\theta),$$

*is satisfied in the notations of the lemma 3.*

*Proof.* As the function $f(x)$ is continuous then it is bounded: there exist a positive constant $M$ such that $|f(x)| \leq M, x \in \Omega$. Taking sufficiently small positive number $\varepsilon$ we can find sylindrical sets $A$ and $A'$ such that $A' \subset B \subset A$ and $\mu_0^*(A \setminus A') \leq \varepsilon$. Then for sufficiently large $T$

$$\frac{1}{T} \int_{\substack{\{t\Lambda\} \in B, \\ 0 \leq t \leq T}} f(\{t\Lambda\}) dt = \frac{1}{T} \int_{\substack{\{t\Lambda\} \in A \\ 0 \leq t \leq T}} f(\{t\Lambda\}) dt + \delta R,$$

where $|\delta| \leq 1$ and

$$|R| \leq M \mu_0^*(A \setminus B) \leq M\varepsilon.$$

Since, $\varepsilon \to 0$ but $M$ is bounded, then passing to the limit as $\varepsilon \to 0$ we get the suitable result as a consequence of the lemma 6.

The statement of a following lemma is an easy consequence of the lemma 7.

**Lemma 8.** *For any function $F(x)$, being continuous in an open set containing the closed set $M \subset \Omega$, the following relation*

$$\lim_{T \to \infty} \frac{1}{T} \int_{\{t\Lambda\} \in M} F(\{t\Lambda\}) dt = \int_M F(\theta) \mu(d\theta)$$

*is satisfied.*

## 4. Proof of the theorem 1.

**Lemma 9.** *Let $L_1(\Omega, \mu_0)$ denote the Lebesgue class of summable functions. Then there exist a positive real function $c = c(\sigma)$ such that*

$$\lim_{T \to \infty} \frac{1}{T} \int_0^T |h(\{t\Lambda\}, \sigma)| \, dt \leq c(\sigma)$$

*for the function $h(\theta, \sigma) = \sum_{k \geq 1} h_k(\theta, \sigma)$ defined above.*

*Proof.* As it has been shown above, the set $A$ of points of divergence of the series

$$h(\theta, \sigma) = \sum_{k \geq 1} h_k(\theta, \sigma)$$

has a zero measure. Let $m$ be some natural number such that $\mu_0(\Omega \setminus U_m) \leq \varepsilon$ for a given small $\varepsilon > 0$. As it was shown in the proof of the lemma 3, in the set $U_m$ the considered series converges uniformly. Therefore, the sum $h(\theta, \sigma)$ of this series is a continuous function in $U_m$, so, it could be found a real number $H = H(\eta)$ such that, everywhere in $U_m$, we have:

$$\left| h(\theta, \sigma) - \sum_{r \leq H} h_r(\theta, \sigma) \right| \leq \eta,$$

for given arbitrarily small positive number $\eta$. Denote by $V$ a closed set $V \subset U_m$ such that $\mu_0(\Omega \setminus V) \leq 1.5\varepsilon$, and by $W$ ($W \supset V$) it's open covering with the condition $\mu_0(W \setminus V) < 0.25\varepsilon$, and also $B = \Omega \setminus V$. Since the set $V$ is closed, then there is a finite family of open balls $(B_r)_{1 \leq r \leq R}$, for which $V \subset \bigcup_{r \leq R} B_r \subset W$ and we can suppose that $W = \bigcup_{r \leq R} B_r$. It is clear that the statement of the lemma 7 is true for the finite union of open balls. So, denoting $\varphi(\theta) = \sum_{r \leq H} h_r(\theta, \sigma)$ we get:

$$\lim_{T \to \infty} \frac{1}{T} \int_{\{t\Lambda\} \in W} \varphi(\{t\Lambda\}) dt = \int_W \varphi(\theta) d\theta.$$

Analogically, there exist an elementary set $W' \subset V$ satisfying the condition $\mu_0(V \setminus W') < 0.25\varepsilon$ that

$$\lim_{T \to \infty} \frac{1}{T} \int_{\{t\Lambda\} \in W'} \varphi(\{t\Lambda\}) dt = \int_{W'} \varphi(\theta) d\theta.$$

So, we get

$$\lim_{T \to \infty} \frac{1}{T} \int_{\{t\Lambda\} \in W'} \varphi(\{t\Lambda\}) dt \leq \lim_{T \to \infty} \frac{1}{T} \int_{\{t\Lambda\} \in V} \varphi(\{t\Lambda\}) dt \leq \lim_{T \to \infty} \frac{1}{T} \int_{\{t\Lambda\} \in W} \varphi(\{t\Lambda\}) dt,$$

from which one deduces using the estimate (6):

$$\left|\lim_{T\to\infty}\frac{1}{T}\int_{\{t\Lambda\}\in V}\varphi(\{t\Lambda\})dt-\int_V\varphi(\theta)d\theta\right|\le\int_{W\setminus W'}|\varphi(\theta)|d\theta\le C_0\sqrt{\varepsilon}\,;\quad(7)$$

here $C_0$ is some constante. Therefore, we have

$$\left|\lim_{T\to\infty}\frac{1}{T}\int_{\{t\Lambda\}\in V}h(\{t\Lambda\},\sigma)dt-\int_V h(\theta,\sigma)d\theta\right|\le 2\eta+C_0\sqrt{\varepsilon}\,.$$

Since the function $h(\theta,\sigma)$ is continuous in the subsets $U_k$ then it is a continuous function in $U$.

Let's consider now a sequence $(g(\sigma\theta))_{\sigma\in\Sigma}$ and define the function $\inf\{g(\sigma\theta)\,|\,\sigma\in\Sigma\}=$
$=\rho(\theta)$ for each considered $\theta$. This function is a measurable function. Define now subsets

$$\Omega_k=\{\theta\in B\,|\,\rho(\theta)\ge 2^k\},$$

where $B$ is an open covering coinciding with a union of open balls $(B_r)_{1\le r\le R}$. Since for each value of $\theta\in\Omega_k$ we have $g(\theta)>2^k$, then from (7) we deduce

$$2^{2k}\mu(\Omega_k)<\int_{\Omega_k}g^2(\theta)\mu(d\theta)\le\left(1+\sum_{r=1}^\infty 2^{-(r-1)\delta}\right)C(\alpha-\delta/2).\quad(11)$$

So we have $\mu(\Omega_k)\le b\cdot 2^{-2k}$ for some positive constant $b$. We can estimate the Lebesgue measure of a curve (which is placed in the closed and finite-symmetrical subset)

$$E_k=\{t\,|\,0\le t\le 1\wedge g(\sigma\{t\Lambda\})\ge 2^k\},$$

by using any finite permutation $\sigma$ satisfying the conditions of the lemma 4 of [13] as a value

$$m(E_k)\le c_0 b 2^{-2k}.\quad(12)$$

Consider now the cases $2^{-k}>\varepsilon$ at points $\theta\in B$ on which $2^{k+1}>\rho(\theta)\ge 2^k$. In this case, to every such point we can put in correspondence some finite permutation, $\sigma$ satisfying the condition

$$g(\sigma\theta)<2^{k+1}$$

in some neighborhood $W(\theta)$ of the point $\theta$. If we define the compact subset $M=\{\theta\in\Omega\,|\,g(\theta)\le 2\varepsilon^{-1}\}$ then it can be overlapped by the union of all open neighborhoods $W(\theta)$ of the points $\theta$. Therefore, there is only finite family of neighborhoods the union of which contains the set $M$. Since, the function $g(\sigma\theta)$ is continuous in the union $(\Omega\setminus B)\cup M$, for each permutation $\sigma$ (this permutation is defined by taken neighborhood $W(\theta)$), then in the set $(\Omega\setminus B)\cup M$ the lemma 7 is applicable:

$$\lim_{T\to\infty}\frac{1}{T}\int_{\{t\Lambda\}\in(\Omega\setminus B)\cup M}g(\{t\Lambda\})dt=\int_{(\Omega\setminus B)\cup M}g(\theta)\mu(d\theta).\quad(12)$$

Consider now the set $B \setminus M$ which is overlapped by the union $A \cup \left(\bigcup E_k\right)$ taken over integers $k$ satisfying the constraints $2^{-k} \leq \varepsilon$ and $A$ is a subset $\Omega$ where the series $\sum_r |g_r(\theta) - g_{r-1}(\theta)|$ diverges.

Let $T_0$ be taken so that, for every $T > T_0$,

$$\left| \frac{1}{T} \int_{\{t\Lambda\} \in (\Omega \setminus B) \cup M} g(\{t\Lambda\}) dt - \int_{(\Omega \setminus B) \cup M} g(\theta) \mu(d\theta) \right| < \eta.$$

Then from (12) it follows that

$$\left| \sum_{n \leq T} \int_{\substack{\{t\Lambda\} \in B \setminus M \\ n \leq t \leq n+1}} g(\{t\Lambda\}) dt \right| < 4c_0 b\varepsilon(T+1).$$

Therefore, for every $T > T_0$

$$\left| \frac{1}{T} \int_{\substack{\{t\Lambda\} \in B \setminus M \\ t \leq T}} g(\{t\Lambda\}) dt \right| < 5c_0 b\varepsilon.$$

Then summarizing all calculations made above, we get the statement 1) of the lemma 7.

To prove the second statement we use the proven fact that for any $T \geq 1$, taking a rectangle $\alpha + \delta \leq \sigma \leq \alpha + 1$, $\alpha = \sigma_m$, $0 \leq t \leq T$, we note satisfiability of the conditions of the lemma 4 for the function $F(s)$. So, this function is an analytical function inside of taken rectangle. This function coincides with given function $f(s)$ in the half plane of absolute convergence (i, e. when $\sigma > \alpha + 1/2$) and, therefore by the principle of analytical continuation $f(s) = F(s)$ in the rectangle. Thus, the relation of the theorem is proven for the case $k = 1$.

Consider now the Dirichlet series $g(s) = \sum_{n=1}^{\infty} b_n n^{-s}$; $b_n = \sum_{n_1 \cdots n_{2k} = n} a_{n_1} \cdots a_{n_{2k}}$. It is easily to see that $g(s) = f^{2k}(s)$. Hence,

$$|b_n|^2 \leq \tau_{2k}(n) \sum_{n_1 \cdots n_{2k} = n} |a_{n_1} \cdots a_{n_{2k}}|^2$$

(here $\tau_k(n)$ expresses the number of representations of $n$ as a product of $k$ natural factors) and we have:

$$\sum_n |b_n|^2 n^{-2\sigma} \leq c_0(\varepsilon) \left( \sum_n |a_n|^2 n^{-2(\sigma-\varepsilon)} \right)^{2k} \leq$$

$$\leq c_0(\varepsilon) C^{2k}(\alpha - (\delta/2)),$$

as $\varepsilon < \delta$, since the following inequality is fair

$$\tau_{2k}(n) \le c_0(\varepsilon) n^\varepsilon$$

with some positive function $c_0(\varepsilon)$ depending only on $\varepsilon$ (see [24, p. 34]). Further, let $r$ be any natural number. Consider the series

$$\sum_{n \in N(r)} b_n n^{-s}.$$

We can write

$$\sum_{n \in N(r)} |b_n| n^{-\sigma} \le \left( \sum_{n \in N(r)} |a_n| n^{-\sigma} \right)^{2k},$$

as the series $\sum_{n \in N(r)} a_n n^{-s}$ is absolutely convergent in the considered rectangle. Applying to the function $g(s) \in J(\sigma_m)$ proved above lemma 9, we get the first statement of the theorem 1, if to accept

$$\lambda_n = \nu_{p_n} = \frac{1}{2\pi} \log p_n,$$

where $p_n$ designates the $n$-th prime number.

The second statement of the theorem 1 is an easy consequence of the first one (see for example [23, p 329]. The proof of the theorem 1 is finished.

### 5. Proof of the theorem 2.

**Lemma 10.** *Let the conditions of the theorem 1 be satisfied for the function $f(s)$, and $f(s_0) = 0$ for some $s_0 = \beta_0 + i\gamma_0$ where $\beta_0 > \sigma_m$. Then, there exist a constant $A$ such that the function $f(s)$ has more than $AT$ zeroes $\beta_0 + i\gamma$ with $|\gamma - \gamma_0| \le T$.*

*Proof.* Let $\delta = (\sigma_m - \beta_0)/4$ and $r < \delta$. Consider a disc $K : |s - s_0| \le r + r/2$ does not containing the zeroes of $f(s)$ different from the $s_0$. Denoting $m_0 = \min_{|s-s_0|=r} |f(s)|$, we apply the lemma 4. Then there exist a natural number $k$ such that

$$\max_{s \in K} |f_k(s) - f(s)| \le 0.01 m_0. \quad (9)$$

Since $f_k(s) = g_k(0; s)$, and the function $g_k(\theta; s)$ is continuous, one can find some cube $C : |\theta_j| \le$
$\le \lambda, j = 1, \ldots, \pi(2^k)$ of the dimension $\pi(2^k)$ such that

$$\max_{s \in K} |f_k(s) - g_k(\theta; s)| \le 0.01 m_0 \quad (10)$$

in the $C$. Consider now the integral

$$I = \frac{1}{2T} \int_{|t| \le T, \{t\Lambda\} \in C} \left[ \iint_K |f(s+it) - f(s)| d\beta d\tau \right] dt$$

where $\{t\Lambda\} = (\{t\log p / 2\pi\})_{p \leq 2^k}$. We have

$$I \leq I_1 + I_2 + I_3 \quad (11)$$

where

$$I_1 = \frac{1}{2T} \int\limits_{|t| \leq T, \{t\Lambda\} \in C} \left[ \iint\limits_K |f(s+it) - f_L(s+it)| d\beta d\tau \right] dt,$$

$$I_2 = \frac{1}{2T} \int\limits_{|t| \leq T, \{t\Lambda\} \in C} \left[ \iint\limits_K |f_L(s+it) - g_k(\{t\Lambda\};s)| d\beta d\tau \right] dt,$$

$$I_3 = \frac{1}{2T} \int\limits_{|t| \leq T, \{t\Lambda\} \in C} \left[ \iint\limits_K |g_k(\{t\Lambda\};s) - f(s)| d\beta d\tau \right] dt,$$

and $L > k$ is a great enough natural number. The first integral we estimate by using of the result of the theorem 1. Since $f(s) \in J = J(\sigma_m)$ then this function has a mean value, so that by the theorem 1 for great enough values of $T$ one has

$$I_1 \leq \frac{1}{2T} \int\limits_{|t| \leq T} \left[ \iint\limits_K |f(s+it) - f_L(s+it)| d\beta d\tau \right] dt \leq$$

$$\leq \iint\limits_K d\beta d\tau \frac{1}{T} \int\limits_{-T}^{T} |f(s+it) - f_L(s+it)| dt \leq 4\pi r^2 \left( \sum_{n \in \mathbf{N} \setminus N(2^L)} |a_n|^2 n^{-2\sigma_0 + 5r/2} + o(1) \right)^{1/2} \leq$$

$$\leq 2^{-\delta L/16} 4\pi r^2 \left( C(\sigma_m + \delta/16) \right)^{1/2}, \quad (12)$$

where $C(\alpha) = \sum |a_n|^2 n^{-2\alpha}$.

For estimating of the second addend of the equality (11) we apply the lemma 6 noting that $f_L(s+it) = g_L(\{t\Lambda\};s)$ (let $\mu = \pi(2^k), \nu = \pi(2^L)$):

$$I_2 = \frac{1}{2T} \int\limits_{|t| \leq T, \{t\Lambda\} \in C} \left[ \iint\limits_K |f_L(s+it) - g_k(\{t\Lambda\};s)| d\beta d\tau \right] dt \leq$$

$$\leq 4\pi r^2 \max_{s \in K} \left[ o(1) + \int\limits_C d\theta_2 \cdots d\theta_{p_\mu} \int\limits_0^1 \cdots \int\limits_0^1 |g_L(\theta;s) - g_k(\theta;s)|^2 d\theta_{p_{\mu+1}} \cdots d\theta_{p_\nu} \right]^{1/2} \leq$$

$$\leq 4\pi r^2 \lambda^{\pi(2^k)} \left( \sum_{n > 2^k} |a_n|^2 n^{-2\sigma_m - 4r} \right)^{1/2}. \quad (13)$$

Consider, at last the integral $I_3$. Applying the lemma 6 and taking into attention the relations (9) and (10), we find:

$$I_3 \leq o(1) + 4\pi r^2 \max_{s \in K} \int\limits_C \left( |f_k(s) - g_k(\theta;s)| + |f(s) - f_k(s)| \right) d\theta_2 \cdots d\theta_{p_\mu} \leq$$

$$\leq 0.002\pi r^2 m_0 \lambda^{\pi(2^k)}. \quad (14)$$

If $k$ taken so large that (note that $\pi(2^k) > 2^k/k$)

$$\left(\sum_{n>2^k}|a_n|^2 n^{-2\sigma_m-4r}\right)^{1/2} \le 0.001\pi r^2 m_0 \lambda^{\pi(2^k)},$$

and $L$ taken as satisfying the condition

$$2^{-\delta L/8} 4(C(\sigma_m+\delta/16))^{1/2} \le 0.001\pi r^2 m_0 \lambda^{\pi(2^k)},$$

we get, combining the relations (9)-(14):

$$I \le 0.01\pi r^2 m_0 \lambda^{\pi(2^k)}. \quad (15)$$

Dissecting the interval $-T \le t \le T$ into more than $2T-2$ unitary subintervals and numbering them we will divide the set of these subintervals into two classes, taking subintervals with odd and even numbers. Then, from (10) one deduces

$$I = \frac{1}{2T}\int_{t\in E',\{t\Lambda\}\in C}\iint_K |f(s+it)-f(s)|d\beta d\tau dt +$$

$$+\frac{1}{2T}\int_{t\in E'',\{t\Lambda\}\in C}\iint_K |f(s+it)-f(s)|d\beta d\tau dt \le 0.01\pi r^2 m_0 \lambda^{\pi(2^k)},$$

where $E'$ and $E''$ imply, correspondingly, the union of subintervals with odd and even numbers. By the lemma 5, one of these sets, for example $E'$, has a measure no less than $(1+o(1))\lambda^{\pi(2^k)}T \ge 0.4\lambda^{\pi(2^k)}T$, if $T$ is sufficiently large. So we have

$$\frac{1}{2T}\int_{t\in E',\{t\Lambda\}\in C}\iint_K |f(s+it)-f(s)|d\beta d\tau dt \le 0.01\pi r^2 m_0 \lambda^{\pi(2^k)}.$$

Denote by $E \subset E'$ subset in which $\iint_K |f(s+it)-f(s)|d\beta d\tau \ge 0.2\pi r^2 m_0$. Then for the measure of this subset we get the bound

$$m(E) \le 0.1\lambda^{\pi(2^k)}T.$$

Then, $m(E'\setminus E) > 0.2\lambda^{\pi(2^k)}T$. Therefore, there exist no less than $[0.2\lambda^{\pi(2^k)}T]-1 = M$ points $-T \le t_1 < \cdots < t_M \le T$, for every of which the following inequality is satisfied:

$$\iint_K |f(s+it_j)-f(s)|d\beta d\tau \le 0.2\pi r^2 m_0 \quad (16)$$

Since the circle $|s-s_0|=r$ can be overlapped by the disks of radiuses $0.5r$ then from (16) one has (see [22, p.345])

$$\max_{|s-s_0|=r}|f(s+it_j)-f(s)| \le 0.8 m_0.$$

Therefore, by the theorem of Rouch'e in every disc $|s+it_j-s_0| \le r, j=1,...,M$ the function have a zero. Since the radius is taken arbitrarily, the proof of the lemma 10 is finished.

For the proof of the theorem we can suppose $a_1 = b_1 = 1$. Really, since $f(s)f(s)^{-1} = 1$ then $a_1 b_1 = 1$ and one can consider the series $a_1^{-1} f(s)$ instead of $f(s)$. Let $0 < X < T$ be any real (unbounded) numbers. Let us to introduce into consideration the series $f(s)\mathrm{M}_X(s)$, where $\mathrm{M}_X(s) = \sum_{n \leq X} b_n n^{-s}$. Every zero of the function $f(s)$ is at the same time a zero of the function $f(s)\mathrm{M}_X(s) \in J(\sigma_m)$. We have

$$f(s)\mathrm{M}_X(s) = 1 + \sum_{n > X} d_n n^{-s}; d_n = \sum_{d|n, d \leq X} b_d a_{n/d}.$$

Further, using known estimations for the number of divisors (see [24]) we get:

$$\sum_{n > X} |d_n|^2 n^{-2\sigma} = \sum_{n > X} \left| \sum_{\delta|d, \delta \leq X} b_\delta a_{d/\delta} \right|^2 n^{-2\sigma} =$$

$$\sum_{n > X} \left| \sum_{\delta|d, \delta > X} b_\delta a_{d/\delta} \right|^2 n^{-2\sigma} \leq \sum_{n > X} \tau(n) n^{-2\sigma} \sum_{\delta|d, \delta > X} |b_\delta|^2 |a_{d/\delta}|^2 \leq$$

$$\leq C_0(\varepsilon) C(\sigma - \varepsilon) \left( \sum_{n > X} |b_n|^2 n^{-2\sigma + 2\varepsilon} \right),$$

where $\sigma > \sigma_m - \varepsilon$ and $C_0(\varepsilon)$ depends only on $\varepsilon$. So, the estimated sum is an infinitesimal as $X \to \infty$.

Let now $s_0$ is a zero placed in the half plane of mean values. By the lemma 10, on the segment $\Re s = \Re s_0, |t - \Im s_0| \leq T$ there are more than $c_0 T$ zeroes, with the distances $\geq 1$ between their ordinates. Let $\rho_1, ..., \rho_M, M \geq c_0 T$ are zeroes placed on the segment $\Re s = \beta_0 > \sigma_m, |t| \leq T$. We have for every zero $\rho$:

$$1 \leq \left| \sum_{n > X} d_n n^{-\rho} \right|^2.$$

Summing over all taken zeroes and applying the lemma 1 of the chapter V of [33] (note that the function $d(f(s)\mathrm{M}_X(s))/dt = -i d(f(s)\mathrm{M}_X(s))/ds \in J(\sigma_m)$ satisfies the conditions of the theorem 1), we find

$$0.5 c_0 T \leq T\, C_0(\varepsilon) C(\beta - \varepsilon) \left( \sum_{n > X} |b_n|^2 n^{-2\sigma + 2\varepsilon} \right).$$

Here $\varepsilon > 0$ is sufficiently small and $X > 1$ great enough. Since the expression parenthesized on the right side is an infinitesimal as $X \to \infty$ then we arrived at a contradiction. The proof of the theorem 2 is finished.

As it said above many of widely used Dirichlet series belong to the class $J$. But the theorem 1 is applicable to the series converging in the half plane of mean values only. When the series has in this half plane singularities, for example the poles, the method used above is applicable for them with some modifications.